\documentclass{svproc}
\usepackage{url}

\usepackage{amsmath,amssymb}
\newcommand{\upPhi}{\mathrm{\Phi}}
\newcommand{\upLambda}{\mathrm{\Lambda}}
\newcommand{\upD}{\mathrm{d}}

\begin{document}
\mainmatter              
\title{Probability Theory in Statistical Physics, Percolation, and Other Random Topics: The Work of C.~Newman}
\titlerunning{The Work of C.~Newman}  
%
\author{Federico Camia\inst{1} \and Daniel L. Stein\inst{2}
}
\authorrunning{Federico Camia and Daniel L. Stein} 
%
\tocauthor{Federico Camia and Daniel L. Stein}
\institute{
	Division of Science, NYU Abu Dhabi, Saadiyat Island, Abu Dhabi\\ and Department of Mathematics, VU University Amsterdam, De Boelelaan 1084A, 1084 HV Amsterdam, The Netherlands,\\
\email{federico.camia@nyu.edu}
\and
	Department of Physics and Courant Institute of Mathematical Sciences,
	New York University, New York, NY 10012 USA\\ 
	and NYU-ECNU Institutes of Physics and Mathematical Sciences at NYU Shanghai, 3663 Zhongshan Road North, Shanghai, 200062, China\\
	and SFI External Professor, Santa Fe Institute, 1399 Hyde Park Rd., Santa Fe, NM 87501 USA,\\
\email{daniel.stein@nyu.edu}
}

\maketitle              

\begin{abstract}
In the introduction to this volume, we discuss some of the highlights of the research career of Chuck Newman. This introduction is divided into two main sections, the first covering Chuck's work in statistical mechanics and 
the second his work in percolation theory, continuum scaling limits, and related topics.
\end{abstract}

\section{Equilibrium and Nonequilibrium Statistical Mechanics}
\label{sec:statmech}

Chuck has devoted a substantial portion of his research career --- and made foundational advances --- to the investigation of statistical mechanical systems, both homogeneous and inhomogeneous, elucidating both their thermal behavior in equilibrium and their nonequilibrium dynamical evolution following a deep quench. Here we briefly review only a few of his most important contributions, beginning with his papers on equilibrium thermodynamics.


\subsection{Thermodynamics of Disordered Systems}
\label{subsec:disordered}

In this section we focus on three areas to which Chuck has contributed heavily, all largely focused on the thermodynamic structure of short-range spin glasses in finite dimensions. These are: general principles of pure state organization (using the metastate approach); multiplicity/non-multiplicity of pure and ground states in realistic and non-realistic models; and presence or absence of a thermodynamic phase transition in sufficiently high dimensions.

{\it General principles of pure state organization}. In the mid-1990's there were (and as of this writing, still are) two leading (and very different) scenarios for the thermodynamic structure of short-range spin glasses: the many-state mean-field replica symmetry breaking~(RSB) scenario put forward by Giorgio Parisi and co-workers~\cite{Parisi79,Parisi83,MPSTV84a,MPSTV84b,MPV87} and the two-state droplet-scaling picture introduced by Macmillan, Bray and Moore, and Fisher and Huse~{\cite{Mac84,BM85,BM87,FH86,FH88b}. Given the analytical intractability of short-range spin glass models, studies relied (and continue to rely) largely on numerical simulations (for a sampling, see~\cite{BY86,CPPS90,RBY90,BCPRPR93,MPR94,MPRR96,CMP97,MPR97,MBD98,PY99a,PY99b,Middleton99,BBDM99,PY00,KM00,DBMB00,BDM00,MP00,MPRRZ00,Middleton00,KPY01,MP01,HYD04,KY05,CGGPV07,KK07,YKM12,Middleton13,BMMMY14,WMK14}). The actual pure state structure that would --- or could --- result from the application of RSB to short-range spin glasses had not been studied in depth; people were mostly looking for numerical evidence (often using spin overlap functions) of the signature features of RSB, in which 1) the thermodynamics is characterized by a mixed thermodynamic state decomposable into a countable infinity of pure states of varying weights within the mixed state; 2) the spin overlap distribution of these pure states is characterized by non-self-averaging (over the coupling realizations) in the thermodynamic limit; and 3) the pure states are ordered hierarchically, in that their relative distances satisfy the ultrametric property~(reviews can be found in~\cite{BY86,RTV86,MPV87,NS03b,SNbook,Read14}). This set of properties taken together at face value was called the ``standard SK picture'' in~\cite{NS97,NSBerlin,NS98,NS03b}.
	
	One of the distinguishing features of Chuck's work with one of us (DLS) on this problem was the use of rigorous mathematical techniques in an area where they had been seldom used before (and perhaps rarer still, where they were used to resolve --- or at least sharpen understanding of --- open physical questions). In the earliest foray on this problem~\cite{NS92}, it was shown that the presence of distinct multiple pure states led to ``chaotic size dependence'': an infinite sequence of finite-volume Gibbs states, generated using boundary conditions chosen independently of the couplings (e.g., the standard periodic, free, or fixed boundary conditions), would not converge in the limit to a single thermodynamic state. However, one could still generate a thermodynamic state by appropriate averaging procedures, and if there were many pure states, this would be a mixture as required by RSB. A surprising conclusion then followed: no matter how the state was constructed (as long as it satisfied basic requirements, such as measurability), the main features of RSB were incompatible with each other.  That is, short-range spin glasses could not support the standard RSB picture in any finite dimension.
	
	This then led to the question, is there {\it any} scenario that is mean-field-like, and that {\it can} be supported in a mathematically consistent fashion in short-range spin glasses? In order to answer this question, the tool of the {\it metastate} was introduced in~\cite{NS96c,NS97,NS98,NS03b}.  (It was soon proved in~\cite{NSBerlin} that this construction was essentially equivalent to an earlier construct introduced by Aizenman and Wehr~\cite{AW90}.) Just as in spin models a thermodynamic state is a probability measure on spin configurations, the metastate is a probability measure on the thermodynamic states themselves. Although it can be used for any thermodynamic system, it is most useful (perhaps even essential) when dealing with systems with many competing thermodynamic states, providing an elegant means of connecting the behavior observed in finite volumes with the (infinite-volume) thermodynamics. (In homogeneous systems, the connection between the two is straightforward, but is anything but straightforward if many pure states exist.) Because of this, its significance and usefulness extend well beyond spin glasses alone (see, for example,~\cite{Kuelske97,Kuelske98}). 
	
	Using the metastate, a ``maximal'' mean-field picture was constructed~\cite{NS96c,NS97,NS03b} that avoided the inconsistencies of the standard RSB model (leading to the moniker ``nonstandard RSB model''). In a series of publications, Newman and Stein (NS) proved that if there are multiple spin glass pure/ground states in finite dimensions, then there must be an {\it uncountable} infinity of them~\cite{NS06b}. In order to arrange them in an RSB-like fashion, there must then be (more formally, the metastate must be supported on) an uncountable set of distinct mixed thermodynamic states, {\it each one of which} is supported on a countable infinity of pure states with varying weights. Roughly speaking, in a given large finite volume, one would observe a single one of these mixed states, which appear with probabilites assigned by the metastate. Ultrametricity would hold among the pure states within a {\it single} mixed thermodynamic state, but not in general between any three pure states chosen arbitrarily from the metastate itself. Moreover, the meaning of non-self-averaging would change, from averaging over different coupling realizations to (roughly speaking) averaging over different volumes for {\it fixed} coupling realization. This picture, based on what remains logically possible based on self-consistency as determined by rigorous arguments, was recently shown by Read~\cite{Read14} to follow from a field-theoretic approach of direct application of RSB techniques to short-range spin glasses.
	
	In very recent work~\cite{ANS15} with Louis-Pierre Arguin, a new set of thermodynamic identities --- connecting pure state weights, correlation functions, and overlaps --- was derived, again using rigorous arguments. It was shown using these relations that pure state weights within a single mixed thermodynamic state consisting of infinitely many pure states (if such a thing exists) are distributed according to a Poisson-Dirichlet process --- exactly as derived much earlier~\cite{MPV85,DT85,Ruelle87}  for the RSB solution to the infinite-range Sherrington-Kirkpatrick~\cite{SK75} spin glass.
	
	So, looking at the larger picture, this body of work used rigorous arguments to greatly narrow down the set of possibilities for the organization of pure states in spin glasses (if in fact there are multiple pure states).  It is important to recognize, as this work emphasizes, that {\it a priori}, there are many possible many-state pictures; having many pure states is not synonymous with RSB. In fact,  the RSB scenario is quite special, having an enormous amount of structure. One of the conclusions of the work briefly summarized above is that {\it if} there are many pure states in short-range spin glasses, then (regardless of dimension) there must be an uncountable number, and {\it if} they appear in nontrivial mixed thermodynamic states (each comprising a countable infinity with given weights), then the {\it only} possibility is the RSB one (as understood within the ``nonstandard'' picture). This is a strong conclusion, and is based on rigorous arguments.
	
	Of course, none of this addresses the actual question of whether the RSB picture (at this point, it's no longer necessary to refer to the ``nonstandard RSB picture'') actually occurs in short-range spin glasses in any finite dimension. In~\cite{NS98,NS03b,SNbook}, it was argued that this is unlikely. The arguments are based on a few rigorous results: the invariance of the metastate with respect to changes of boundary conditions~\cite{NS98}, and important differences in the behavior of edge disorder chaos in infinite-range and short-range models, as pointed out by Chatterjee~\cite{Chatterjee09}. (We note that the argument claiming that metastate invariance is incompatible with RSB has been criticized in~\cite{Read14}.) Perhaps more compelling is a nonrigorous argument based on pathologies arising from coupling strengths in infinite-range spin glasses scaling to zero with the number of spins, coupled with the randomness of coupling signs~\cite{NS03c,SNbook}. However, the issue remains unresolved as of now. 
	
	While the body of work described above has led to the ruling out of an enormous number of possible scenarios, it has also led to the introduction of new ones. A very natural and (so far) viable picture, introduced and investigated in~\cite{NS96c,NS97,NS98,NS03b} and called the ``chaotic pairs'' picture, arises naturally from the metastate approach. Chaotic pairs is a scenario that, like RSB, has uncountably many pure states in the (periodic boundary condition, for specificity) metastate, but they are organized much more simply than in RSB, with a simple, self-averaging overlap structure that is identical (in any finite volume) to that of droplet-scaling. In chaotic pairs, thermodynamic states are a {\it trivial} mixture of pure states, each consisting of a single pure state and its global flip, each with weight 1/2 in the mixed thermodynamic state. However, there are uncountably many of these thermodynamic states, each one comprising a  pure state pair different from all the others. 
	
	For completeness, we should mention a fourth picture, introduced independently by Palassini and Young~\cite{PY00} and Krzakala and Martin~\cite{KM00} and sometimes referred to as the ``TNT'' picture. Here excitations above the ground state have boundaries with zero~density in the edge set (like droplet-scaling, and unlike RSB and chaotic pairs) but whose energies do not scale as some increasing function of the volume (like RSB, and unlike droplet-scaling and chaotic pairs). While this last picture does not specify how many pure/ground states there are, it was argued in~\cite{NSregcong01} that it is most naturally consistent with two-state pictures. If so, that leaves us with two potential many-state pictures for the short-range spin glass phase (RSB, chaotic pairs) and two two-state pictures (droplet-scaling, TNT). Which of these holds for short-range spin glasses, and in which dimensions, remains a fundamental open problem in the mathematics and physics of the statistical mechanics of disordered systems.

	{\it Multiplicity of pure and ground states in short-range spin glasses.}  The work described in the previous section clarified and restricted the ways in which pure states could
	be {\it organized} in short-range spin glasses; it did not address the actual {\it number} of pure or ground states of the spin glass phase. It should be noted that --- from a logical standpoint at least ---
	this question cannot be considered without answering several deeper ones: namely, is there a thermodynamical phase transition at nonzero temperature in {\it any} finite dimension, and if so, does the low-temperature phase break spin-flip symmetry (so that pure states come in spin-reversed pairs, as described in the previous section)? The first of these questions will be addressed in the following section; for now we will assume, in accordance with a good deal of numerical evidence~\cite{BY86,GSNLAI91,PC99,Ballesteros00,HPV08,Janus13} that above some lower critical dimension there is a phase transition at a dimension-dependent positive temperature to a spin glass phase with broken spin-flip symmetry (for those readers interested in more offbeat possibilities, see~\cite{NS06b}). For specificity we restrict the discussion in this section to spin glasses described by the nearest-neighbor Edwards-Anderson model~\cite{EA75} on the Euclidean lattice $\mathbb{Z}^d$ and with i.i.d. couplings taken from the normal (Gaussian) distribution with mean~zero and variance~one.
	
	The main question, if pure states come in pairs, is whether there is only a single pair, as in droplet-scaling and (probably) TNT, or else uncountably many, as in RSB and chaotic pairs. (We will always refer in this section to infinite-volume pure or ground states. Definitions and constructions can be found in~\cite{NS03b}.})  Aren't there other possibilities? Probably not: a countable number, either finite or infinite, was rigorously ruled out for RSB as noted earlier, and there exists strong evidence that it cannot occur for chaotic pairs either. The number of pure/ground states, of course, could in principle be dimension-dependent, or even temperature-dependent for fixed dimension. But we're a long way at this point from addressing these questions. (We confine our discussion to the so-called {\it incongruent} ground/pure states~\cite{FH87,HF87}, which differ from each other by positive-density interfaces, and are generated by coupling-independent boundary conditions~\cite{NSregcong01,NS03b}.)

Until recently the only dimensions in which the answer was fully known were one dimension (a single pure state at all positive temperatures and a single pair of ground states at zero temperature) and (strictly) infinite dimensions (where it is assumed that the state space structure is given by RSB). For the moment, let's assume the simplest (though not the most interesting) case that in each dimension there is a single transition temperature, above which there is a single pure state and below which there is a spin glass phase where the pure state cardinality is independent of temperature and equal to the ground state cardinality. This appears to be what most workers in the field believe. Many who fall into the RSB camp are inclined toward the possibility that at the lower critical dimension (possibly $d=3$, probably no larger than $d=4$) and above, the low-temperature phase is fully described by RSB. Others have argued that there is also an upper critical dimension at $d=6$~\cite{Aspelmeier16}, below which there is a single pair of pure/ground states and above which infinitely many. Still others~\cite{FH87} conjecture that there is only a single pair of pure/ground states in all finite dimensions (above the lower critical dimension), and the infinite-dimensional limit is singular in the sense of the structure of the spin glass phase.

All of this remains conjecture. Numerical experiments mostly (but not universally) agree that in $2D$ there is only a single pair of ground states (it is believed that $T_c=0$ in $2D$), but in $3D$ and $4D$ different groups have arrived at different conclusions based on their numerical studies. The only rigorous results for any dimension between one and infinity appear in~\cite{NS2D00,NS2D01,ADNS10}, all of which deal with two~dimensions. In~\cite{NS2D00,NS2D01}, NS proved that in $2D$ any two distinct ground states can differ by only a {\it single}, positive~density interface, providing evidence that there is only a single pair of ground states. Further evidence for this conclusion is provided in~\cite{ADNS10}, where Arguin, Michael Damron, and NS proved essentially that the {\it half}-plane with a free boundary condition along its edge has only a single pair of ground states.

In a recent set of papers, Arguin, NS, and Janek Wehr attacked this problem from a different, dimension-independent, direction, by proving lower bounds on free energy fluctuations between pure states at positive temperature~\cite{ANSW14} and energy fluctuations between ground states at zero temperature~\cite{ANSW16}.  An upper bound has existed for many years, where it was shown that in finite volumes free energy fluctuations scale no faster than the square root of the surface area of the volume under consideration~\cite{AFunpub,NSunpub,Stein16article}.
Lower bounds have already been obtained for several cases, where it was shown that free energy fluctuations scale linearly with the square root of the volume under consideration.  A (hopefully small) difficulty is that the quantity used to obtain the lower bound is slightly different from that used for obtaining the upper bound, so a corresponding upper bound (which is expected to behave similarly to the known upper bound just described) needs to be found for the quantity examined in~\cite{ANSW14,ANSW16}. If it can, then for any cases where these results can be shown to apply, there must be no more than a single pair of pure/ground states. This work is in progress at the time of this writing.

While progress has been slow in determining numbers of pure/ground states in realistic spin glass models, there are other interesting, but unrealistic, models which can provide some interesting illumination. A very useful such model is the so-called {\it highly disordered} model~\cite{NS94,BCM94,NS96a}, in which the Hamiltonian is Edwards-Anderson and couplings are independent random variables, but chosen from a volume-dependent distribution. The idea is that the distribution of coupling magnitudes depends on volume in such a way that, in sufficiently large volumes, each coupling magnitude is more than twice as large as the next smaller one, and no more than half as large as the next larger one. One interesting feature of this model is that the true ground state in any volume can be found by the greedy algorithm (in fact, this was how the model was originally arrived at). Moreover, the algorithm used to find the ground state can be mapped onto invasion percolation.  Consequently, the number of ground states can be found in any dimension: there is only a single pair below $d=6$ and uncountably many above. (The original papers~\cite{NS94,NS96a} used heuristic arguments to propose $d=8$ as the crossover dimension, but a more careful analysis in~\cite{RJ} indicated the crossover dimension to be $d=6$.)

Interestingly, the ground state multiplicity in this model is the same for both the spin glass and the random ferromagnet. The somewhat subtle distinction between the two comes from features of chaotic size dependence for the spin glass, and its absence in the ferromagnet. In any case, the highly disordered model remains a useful testing ground for new ideas.

{\it Equilibrium phase transitions.}  Not much of the
discussion so far is relevant to spin glasses if there is no equilibrium
phase transition to a spin glass phase above some lower critical dimension
$d_c$: at positive temperature there would simply be a unique~Gibbs
state. (Of course, even in that case the question of multiplicity of ground
states at $T=0$ would still exist.) Even with a phase transition, it could
still be the case that there exists a unique Gibbs state above and below
the transition temperature $T_c(d)$; the state would simply be qualitatively different in
these two temperature regimes. Or there could be multiple transitions
within a fixed dimension~\cite{NS06b}. And even in the
absence of a phase transition there could be some sort of dynamical
transition with interesting features.

However, a great deal of mostly numerical work~\cite{BY86,GSNLAI91,PC99,Ballesteros00,HPV08,Janus13} points toward a more
conventional scenario in which in dimension $d>d_c$, there is a unique transition temperature
$T_c(d)$, below which spin-flip
symmetry is broken (equivalently, the Edwards-Anderson order
parameter~\cite{EA75} $q_{EA}\ne 0$). Moreover, almost all (currently
relevant) theoretical work addressing the low-temperature spin glass phase,
whether using evidence from the infinite-range spin glass (i.e., the RSB
scenario), or resulting from scaling arguments (droplet-scaling), or
suggested from the structure of the metastate (chaotic pairs) or following
numerical simulations (TNT) begin by assuming an equilibrium spin glass
transition breaking spin-flip symmetry. So it's natural to ask whether such
a transition can be proved to exist (in {\it any} finite dimension), and
if so, what can one say about its structure?

Substantial progress toward this end was made by Jon Machta and NS in~\cite{MNS08,MNS09a,MNS09b}, where
percolation-theoretic methods were used to uncover what is likely to be an
underlying geometric structure for spin glass phase transitions. These
random graph methods, in particular the Fortuin-Kastelyn (FK) random
cluster (RC) representation~\cite{KF69,FK72}, provide a set of useful tools
for studying phase transitions (more specifically, the presence of multiple
Gibbs states arising from broken spin flip symmetry) in discrete spin
models. These representations map Ising and Potts models onto a type of
percolation problem, thereby allowing spin correlation functions to be
expressed through the geometrical properties of associated random
graphs. FK and related models are probably best known in the physics
literature for providing the basis for powerful Monte Carlo methods, in
particular the Swendsen-Wang algorithm for studying phase
transitions~\cite{SW87}. At the same time they have proved important in
obtaining rigorous results on phase transitions in discrete-spin
ferromagnetic (including inhomogeneous and randomly diluted) models. Because of complications due to frustration, however,
graphical representations have so far played a less important role in the
study of spin glasses.

Some readers will already be familiar with the Fortuin-Kasteleyn random
cluster representation. For those who are not and would like to read on, a
brief one-page summary can be found in~\cite{MNS09b}. (For those who have both
the interest and the time, the original references~\cite{KF69,FK72} are the
best place to go. Also see the article by Newman~\cite{Newman94} and/or the lengthy and detailed review by Grimmett~\cite{Grimmett06}, as well as papers by Camia, Jianping Jiang, and Chuck containing recent applications of FK~methods~\cite{CJN17a,CJN18+}.) In what follows, we assume a familiarity with the basic ideas of the FK random cluster representation for discrete spin models.  
For our purposes, it's sufficient to note that in the case of a ferromagnet, general theorems~\cite{BK91} ensure that when percolation of an FK cluster occurs, the percolating cluster is unique.  
We shall informally refer to
such percolation as ``FK percolation''. It can then be shown that FK percolation in the ferromagnet corresponds to the presence of multiple Gibbs states (in the Ising ferromagnet, 
magnetization up and magnetization down) with broken spin-flip symmetry, and moreover that the onset of percolation occurs at the ferromagnetic critical temperature. 

However, the situation is greatly complicated when couplings of both signs are allowed, as in the spin glass. In this case, FK percolation does {\it not} appear to be a sufficient condition for 
multiple Gibbs states (although it's undoubtedly necessary), and the numerical onset of FK percolation does not coincide with (what is believed to be) the transition to a spin glass state with broken spin-flip symmetry~\cite{SW87}.
The essential problem is the following. If there is broken spin-flip symmetry, then in a finite-volume Gibbs state one should be able to change the orientation of the spin at the origin by changing the orientations of fixed spins at the boundary; or equivalently, the use of (at least some) fixed-spin boundary conditions should lead to a nonzero thermally averaged magnetization at the origin. As an example, choose the boundary condition where all boundary spins are fixed to be $+1$.
At positive temperature, there are many possible realizations of the FK percolating cluster. But at least {\it a priori}, some of these will connect the origin to the boundary in a way that will lead to a positive magnetization at the origin, and others to a negative magnetization, and these could cancel leading to a net zero thermal average of the magnetization. Of course, we don't know for sure that this happens, but so far no one has been able to rule it out. Supporting evidence comes from the numerical studies described above, in which single FK percolation occurs well above what is believed to be the spin glass transition temperature, and proofs that FK percolation occurs~\cite{GKN92} in $2D$ models where no broken spin-flip symmetry is expected. An interesting speculation, which remains to be studied numerically, is that FK percolation might indeed indicate a phase transition, but not one where spin-flip symmetry is broken. Rather, it could be the case that there is a change in behavior of two-point correlation functions with distance, but the EA~order parameter remains zero~\cite{NS06b}.

So at first glance it would appear that FK percolation is the wrong tool to use in analyzing spin glass phase transitions. However, it turns out that FK percolation is indeed relevant, but that there is a considerably more complicated signature for the spin glass transition in both short-range and infinite-range models, as elucidated in~\cite{MNS08,MNS09a,MNS09b}. These papers pointed out that understanding the percolation signature for the spin glass phase requires two ingredients beyond what is needed for ferromagnets. The first is the need to consider percolation within a two-replica representation, and the second is that spin glass ordering corresponds to a more subtle percolation phenomenon than simply the appearance of a percolating cluster ---  it involves a {\it pair} of percolating clusters, first proposed in~\cite{NS06b}. 

The ferromagnetic phase transition corresponds to percolation of a unique infinite cluster in the FK representation, which will hereafter be referred to as ``single'' FK percolation. As noted above, single FK percolation is insufficient --- so far as is currently known --- to indicate the presence of multiple Gibbs states in spin glasses. However, if one switches to a two-replica formalism, one can study what might be called ``double FK percolation'', meaning the following: take two independent FK realizations and consider ``doubly occupied'' bonds --- i.e., bonds occupied in both representations. If these percolate (as usual, with probability one), then one has double FK percolation. (In fact, there are two-replica representations other than FK than can be used equally well, as described in~\cite{MNS08}. For the sake of brevity we focus here only on double FK percolation.)

Double FK percolation is much more difficult to study than single FK percolation, and a number of general theorems valid for the latter fail for the former. In particular, uniqueness of a percolating double cluster is no longer guaranteed. And not only is it not guaranteed, it doesn't happen. What was found instead was that the spin glass transition corresponds to the breaking of indistinguishability between {\it two} percolating networks of doubly FK-occupied bonds ---  in particular, by their having a nonzero difference in densities.

It's worth describing this in just a little more detail, starting with short-range spin glasses (in what follows, we note that results were obtained numerically for the EA~spin glass in three dimensions and rigorously for the infinite-range SK model). For the EA~model in $3D$ (and presumably in higher dimensions as well), there is a series of transitions as the temperature is lowered from infinity. At a temperature $T_{c2}$ well above the putative spin glass transition temperature, there appear two doubly-infinite FK clusters of equal density, and presumably macroscopically indistinguishable. In terms of using boundary conditions to observe macroscopically distinguishable Gibbs states, this situation is no better than that of single FK percolation. Below a lower temperature $T_{c1}$, which is (within numerical error) equal to the spin glass transition temperature $T_c$ observed using other means, the clusters separate in density --- one grows and the other diminishes, leading to the presence of a macroscopic observable whose sign can be reversed using a flip of boundary conditions (and therefore corresponds to broken spin-flip symmetry).

The SK model behaves differently: above $T_c$, there is no double percolation at all, while below $T_c$ two doubly infinite clusters of unequal density appear.

For more details we refer the reader to~\cite{MNS08,MNS09a,MNS09b}. This work will hopefully prove useful in understanding better the differences between the nature of the phase transition in ferromagnets and in spin glasses. More importantly, the hope is that for short-range models, this work will be a significant step toward finally developing a rigorous proof for spin glass ordering and will eventually lead to a clean analysis of the differences between short- and infinite-range spin glass ordering. Finally, it can help to explain why there is no spin glass transition leading to broken spin-flip symmetry on simple planar lattices: two dimensions does not generally provide enough ``room'' for two disjoint infinite clusters to percolate. However, a system that is infinite in extent in two dimensions but finite in the third might be able to support two percolating clusters, with unequal densities at low temperature. These intriguing possibilities remain unexplored at the time of this writing.

\subsection{Nonequilibrium Dynamics of Discrete Spin Systems}
\label{subsec:dynamics}

A second broad area of Chuck's research interests centers on the dynamical evolution of Ising and Potts-like systems, both ordered and disordered, in conditions far from thermodynamic equilibrium. His work in this general area addresses a wide variety of problems focusing on different aspects of nonequilibrium dynamics, but this section will touch on only two: the mostly unexplored relation between nonequilibrium dynamics and equilbrium thermodynamics, and the problem of predictability (or more colloquially, ``nature vs.~nurture'') in spin systems following a deep quench. For brevity's sake, we unfortunately omit other areas of Chuck's notable contributions to topics as diverse as persistence~\cite{NS99b}, aging~\cite{FIN99,FINS01,FIN02}, broken ergodicity~\cite{NS95}, biological evolution~\cite{KN85,NCK85}, and food webs~\cite{NC85,CNB85}; we refer those interested to the references cited.

{\it Nonequilibrium dynamics and equilbrium thermodynamics.}  It is often noted in the literature that pure state multiplicity, as a purely equilibrium property of thermodynamic systems, plays no role in its dynamics, because (under conventional dynamics, such as one-spin-flip, to which we adhere throughout) a system in a pure state remains in that state forever. A natural conclusion to draw is that a system's pure state multiplicity and its dynamical behavior are largely distinct, and information about one says little about the other. However, a series of rigorous arguments in~\cite{NS99a} showed a surprisingly deep connection between the two under what initially would seem the most unpromising conditions: a system undergoing a deep quench and thereafter evolving under conditions far from equilibrium.

The arguments are long and technical and won't be repeated here. The essence of the main result is the following.  Consider a system at a temperature where, in equilibrium,
there are (or are presumed to be) two or more pure
states.  Examples include the Ising ferromagnet and the EA~spin glass above the critical dimension, among many others. Consider now a deep quench from high (well above $T_c$) to low (well below $T_c$) temperature.
It was proved in~\cite{NS99a} that as time $t \to \infty$, although the system is usually
in some pure state locally (i.e., within a fixed volume), then
either a) it {\it never} settles permanently (within that volume)
into a single pure state, or b) it does but then the pure state depends on {\it both} 
the initial spin configuration {\it and}  the realization of the stochastic dynamics.

It was further proved that the former
case holds for deeply quenched $2D$ ferromagnets.
But the latter case is particularly interesting, because it was shown in~\cite{NS99a} that it can occur only if there exists an uncountable number 
of pure states with almost every pair having zero overlap.
It was further shown that in both cases, almost {\it no} initial spin configuration is in the
basin of attraction of a single pure state. That is, after a deep quench the resulting configuration
space is almost all boundary; equivalently,  the union of the 
basins of attraction of all pure states forms a set of measure
zero in configuration space. 

So in principle the nonequilibrium dynamics following a deep quench provides information about the multiplicity of pure states
and vice-versa. However, this relation is likely to be more useful for future mathematical analysis than for numerical tests, which 
are confined to relatively small systems. Nevertheless, there may be experimental consequences as yet unexplored.
Unlike the ferromagnet, it is not clear that spin glasses are easily prepared to lie within a single pure state, even under conditions as
close to equilibrium as current technology allows. In particular, because of the possibility of chaotic temperature
dependence \cite{BM87,FH88b}, the conclusions of~\cite{NS99a} could well apply to laboratory spin glasses prepared under conditions where small temperature
changes are made slowly.  Experimentally observed slow relaxation
and long equilibration times in spin glasses could then be a consequence of small (relative to the system) domain size
and slow (possibly due to pinning) motion of domain walls (a conclusion earlier reached by Fisher and Huse~\cite{FH88a} using different considerations).

{\it Nature vs.~nurture: Predictability in discrete spin dynamics.}   Although the core issues raised in~\cite{NS99a} remain open, the paper led to
an unanticipated research direction that opened an entirely new area of inquiry: given a typical initial configuration, which then
evolves under a specified dynamics, how much can one predict about the
state of the system at later times?  Chuck's papers with collaborators have colloquially referred to
this as a ``nature~vs.~nurture'' problem, with ``nature"
representing the influence of the initial configuration (and disorder realization, if relevant) and ``nurture''
representing the influence of the random dynamics.

First one must determine under what conditions a system settles down~\cite{NS99a,NNS00}, in the sense of {\it local equilibration}: 
do domain walls cease to sweep across a fixed region after a (size-dependent) finite time? Local equilibration occurs in many
quenched systems, for example, any disordered Ising model with continuous couplings~\cite{NNS00}, but has also been shown {\it
	not} to occur for others, such as the $2D$~homogeneous Ising
ferromagnet~\cite{NNS00}. When local {\it non}equilibration (LNE) occurs, one can still ask
whether the {\it dynamically averaged} configuration has a limiting distribution.
If so, this implies a complete lack of predictability, while the absence of a distributional limit implies that some amount of predictability
remains~\cite{NS99a,NNS00}. In~\cite{NS99a}, systems displaying LNE but having a limiting dynamically averaged distribution were said to exhibit ``weak LNE'', while those
with no distributional limit were said to exhibit ``chaotic time dependence''.

So the distinction between these two types of local nonequilibration is important for the question of predictability. While the presence of~LNE has been rigorously established for
several systems~\cite{NS99a,NNS00}, determining which type of LNE occurs is considerably more difficult (when the initial state is chosen randomly, as discussed above). However, Renato~Fontes, Marco~Isopi,  and Chuck were
able to prove the presence of chaotic time dependence in the one-dimensional voter model with random rates at zero temperature (equivalently, the $1D$ zero-temperature stochastic Ising~model) if the disorder distribution is heavy-tailed~\cite{FIN99}. (This is a $1D$ version of Jean-Phillipe~Bouchaud's trap model~\cite{aging6}). Otherwise, they showed that
chaotic time dependence is absent. The latter conclusion was also shown to hold for the voter model with random rates in dimension greater than two. In further work~\cite{FIN02}, Fontes~{\it et al.} extended the results of~\cite{FIN99} 
by analyzing in detail the space-time scaling limit of the $1D$ voter model with random rates as a singular diffusion in a random
environment. The limit object, now know as FIN (Fontes-Isopi-Newman) diffusion, has some striking scaling properties. Much subsequent work, in particular by Gerard Ben~Arous, Jiri Cerny, and
collaborators~\cite{BC06} has studied scaling limits beyond one dimension, which are rather different kinds of objects than the FIN diffusions.

Returning to the problem of nature vs.~nurture, the questions raised by these considerations have led to a rich area of research that has brought in a variety of collaborators and students, and
that touches on a number of related dynamical problems of interest, including phase-ordering kinetics, persistence, damage spreading, and aging.
The history of these efforts, along with recent advances, is described in detail in a separate contribution to this volume by one of us (DLS), to which we refer
the interested reader.

\subsection{Lee-Yang zeros and the Riemann Hypothesis.} 
\label{subsec:LY}

In 1950 T.D.~Lee and C.N.~Yang published a pair of papers~\cite{LY1,LY2} that pioneered a new way of understanding phase transitions. They considered the ferromagnetic Ising~model in an external magnetic field, and showed that the zeros of its partition function, as a function of the external magnetic field, all lie on the imaginary axis. Since that time it's been shown that a number of other statistical~mechanical systems obey Lee-Yang type theorems of their own, and Chuck played a substantial part in these efforts. Given that Lee-Yang type theorems provide not only an understanding of the properties of the phase transition in a given system (e.g., existence of a mass gap~\cite{PL74,GRS75}), but can also lead to useful correlation inequalities~\cite{Dunlop79,Newman75}, such results are important in our overall understanding of phase transitions in condensed matter. For a general review and discussion of applications of Lee-Yang type theorems, see~\cite{FR12}. 

Among Chuck's contributions to this area was to extend the Lee-Yang theorem to ferromagnetic models with a very general class of single-spin distributions~\cite{Newman74}, to the classical XY~model~\cite{DN75} (further generalized in~\cite{LS81}; see also~\cite{SF71}), and~\cite{FS81,NW17a,NW17b} to Villain~models~\cite{Villain75}, which are closely related to XY~models. In~\cite{NW17b} it was further shown that complex Gaussian multiplicative chaos in general does {\it not} have the Lee-Yang property. These results have some interesting implications. In particular, the ``spin wave conjecture''~\cite{Dyson56,MW66} asserts that, below a critical temperature, the angular variable $\theta$ of the XY (and Villain) model at large scales behaves like a Gaussian free field (modulo $2\pi$), suggesting in turn that the spin variables in these models could behave like a version of complex Gaussian multiplicative chaos. The spin wave conjecture might lead one to expect that complex Gaussian multiplicative chaos would display a Lee-Yang property; but the paper of Newman and Wu rules that out, at least for a range of inverse temperature~$\beta$.

In a different direction, Chuck has long been interested in the possible connection between the Lee-Yang theorem and the famous Riemann Hypothesis. Recall that the Riemann Hypothesis states that the nontrivial zeroes of the Riemann zeta function $\zeta (s)=\sum_{n=1}^\infty\frac{1}{n^s}$ all have real part~$1/2$. While on the surface it sounds less than dramatic, many consider it to be the most important unsolved problem in mathematics. Its proof (or disproof) has important implications for the distribution of prime numbers, the behavior of various functions in number theory and combinatorics, eigenvalue distributions of random matrices, quantum chaos, and problems in many other areas. The similarity between the distribution of Lee-Yang zeros and those of the Riemann zeta function have long led to speculation that one path to a proof of the Riemann Hypothesis lies through the Lee-Yang theorem (for a review, see~\cite{Knauf99}). To date, however, attempts to carry out this program have been unsuccessful (otherwise, you surely would have heard).

However, Chuck made an important contribution to this problem in a slightly different and very interesting direction. In 1950 Nicolaas De Bruijn showed that a certain function $H(\lambda,z)$ (its precise form is unimportant for the purposes of this discussion) has only real zeros for $\lambda\ge 1/2$~\cite{DeB50}. In addition, if $\lambda$ is such that $H(\lambda,z)$ has only real zeros, then for all $\lambda'>\lambda$, $H(\lambda',z)$ also has all real zeros. The Riemann~Hypothesis is equivalent to $H(0,z)$ having only real zeros.  Now if $H$ has only real zeros for {\it all}  real $\lambda$, then the Riemann~Hypothesis would follow. This was a strategy attempted by Polya, among others. However, Chuck proved in~\cite{Newman76} that there exist real $\lambda$ for which $H$ has a nonreal zero. This led to what is now known as the De Bruijn-Newman constant, usually denoted by $\upLambda$. It is defined to be the value such that if $\lambda\ge \upLambda$, $H$ has only real zeros, while if $\lambda<\upLambda$, $H$ has a nonreal zero.

The Riemann Hypothesis is true if and only if the De Bruijn-Newman constant $\upLambda\le 0$. But in~\cite{Newman76}, Chuck conjectured that $\upLambda\ge 0$. If true, this immediately implies that the Riemann Hypothesis is true if and only if $\upLambda=0$. Computer-aided rigorous calculations of a lower bound for the De Bruijn-Newman constant have been made over the years. Until 2018, the best lower bound was $-1.1\times 10^{-11}$~\cite{SGD11},  but recently Brad Rodgers and Terry Tao~\cite{RT18} posted a proof verifying Chuck's conjecture that $\upLambda \ge 0$. 

There also exist upper bounds for $\Lambda$. The earliest was de~Bruijn's $\Lambda\le 1/2$~\cite{DeB50} in his original paper; an improvement to a strict inequality
$\Lambda<1/2$ was made in 2009 by Ki, Kim, and Lee~\cite{KKL09}.  An improved upper bound of~0.22 was very recently obtained by Tao and collaborators, as posted in~\cite{polymath}. In addition, a survey article by Chuck and Wei Wu on de~Bruijn-Newman~type constants has just been published~\cite{NW19}.

While work is ongoing, we'll give the last word to Chuck. In~\cite{Newman76}, Chuck noted that his conjecture that $\upLambda\ge 0$ put a quantitative scaffolding behind the dictum that the ``Riemann hypothesis, if true, is only barely so.''

This introduction has barely described the extent and importance of Chuck's many contributions to mathematical physics and related areas. Our hope is that the reader will go on to take a look at some of Chuck's many papers, especially in areas that for lack of space, and the finite lifetimes of the authors, this introduction was not able to discuss.

In the following sections we turn to a different aspect of Chuck's work, involving seminal contributions to percolation theory, scaling limits, SLE, and related topics.

\section{Percolation Theory and Continuum Scaling Limits}
\label{sec:percolation-scaling}

In this section, we briefly review some of the other fundamental advances, more probabilistic in nature, made by Chuck in his long and productive career. We will only briefly mention the Brownian Web and then discuss Chuck's contributions to percolation theory and to the study of continuum scaling limits. Like the results surveyed in the previous section, this side of Chuck's work was mainly motivated by questions arising in statistical physics, and sometimes quantum field theory, or a combination of the two. The breadth of Chuck's interests and his agile versatility as a mathematician and mathematical physicist are evident in the collection of problems discussed here and in the previous section.

\subsection{The Brownian Web}
\label{subsec:web}

{\it Construction and relation to disordered systems}.
One of Chuck's strengths is his ability to move between fields and exploit fruitful connections. This is evident in his influential work on the Brownian web, which has its roots in the analysis of nonequilibrium dynamics in one-dimensional disordered systems, discussed in the previous section. Roughly speaking, the Brownian web is the scaling limit of the space-time graphical representation of an infinite collection of coalescing random walks. It originated from Arratia's Ph.D. thesis~\cite{Arratia79}, where it is shown that a collection of coalescing random walks on $\mathbb Z$ starting from every vertex of $\mathbb Z$ converges to a collection of Brownian motions on $\mathbb R$ starting from every point of the real line at time zero. Arratia's attempt~\cite{Arratia81} to construct a process corresponding to a collection of coalescing Brownian motions starting at every point of the real line at every time $t \geq 0$ was never completed. However, building on Arratia's work, the process was constructed years later by Balint T\'oth and Wendelin Werner \cite{TW98}, who discovered a surprising connection between Arratia's one-dimensional coalescing Brownian motions and a two-dimensional random process repelled by its own local time profile, which they called {\it true self-repelling motion}. A few more years later, Luiz Renato Fontes, Marco Isopi, Chuck and DLS realized that the same system of coalescing Brownian motions also arises from the scaling limits of one-dimensional spin systems \cite{FINS01}.
The first three authors, together with Ravishankar, introduced~\cite{FINR02,FINR04} a topology such that the system of coalescing Brownian motions starting from every space-time point can be realized as a random variable in a separable metric space, and they named this random variable the {\it Brownian web}. This remarkable object emerges in contexts as diverse as hydrology and the zero-temperature dynamics of the Ising model. The work of Chuck and coauthors on the Brownian web has inspired others to explore its properties and propose extensions, spurring a wealth of interesting papers. We will not discuss the topic further in this introduction; instead we refer the interested reader to~\cite{SS08,SSS17} for a survey of results and extensions.

\subsection{Percolation}
\label{subsec:percolation}

Percolation as a mathematical theory was introduced by Broadbent and Hammersley~\cite{Broadbent54,BH57}
to model the spread of a ``fluid'' through a random ``medium.'' Broadbent and Hammersley interpreted the terms fluid and medium broadly, having in mind situations such as a solute diffusing though a solvent, electrons moving through an atomic lattice, molecules penetrating a porous solid, or disease infecting a community. They were interested in situations where the randomness is associated with the medium rather than with the fluid. To mimic the randomness of the medium, one can think of the latter as a system of channels some of which are randomly declared closed to the passage of the fluid. This can be modeled by a $d$-dimensional cubic lattice, seen as an infinite graph, where the edges between nearest-neighbor vertices are independently declared \emph{open} (to the passage or the fluid) with probability $p$ or \emph{closed} with probability $1-p$. Models of this type are called \emph{bond percolation} models, and can be defined on any graph. Many other variants have been studied, attracting the interest of both mathematicians and physicists. One version in particular, site percolation on the triangular lattice, will be discussed in more detail in Section~\ref{subsec:scaling}. In \emph{site percolation}, the vertices of a graph rather than the edges are declared open or closed.

Mathematicians are interested in percolation because of its deceptive simplicity which hides difficult and elegant results. From the point of view of physicists, percolation is one of the simplest statistical mechanical models undergoing a continuous phase transition as the value of the parameter $p$ is varied, with all the standard features typical of critical phenomena (scaling laws, conformal invariance, universality). On the applied side, percolation has been used to model the spread of a disease, a fire or a rumor, the displacement of oil by water, the behavior of random electrical circuits, and more recently the connectivity properties of communication networks.

{\it Existence and uniqueness of infinite clusters}.
One of the most interesting aspects of percolation, and a major reason for its popularity, is its phase transition, which is of a purely ``geometric'' nature. Defining open clusters to be the sets of vertices connected to each other by a path of open edges, one can ask whether there exists an infinite open cluster at a given value $p$ of the density of open edges. This amounts to asking whether there is an open path from the origin of the lattice (any deterministic vertex) to infinity with strictly positive probability. Indeed, if the answer to this question is positive, then translation invariance and an application of Kolmogorov's zero-one law imply the existence of an infinite open cluster with probability one somewhere in the system. If the answer is negative, then with probability one, no infinite open cluster exists. This justifies the introduction of the \emph{percolation function} $\theta(p)$ defined as the probability that the origin is connected to infinity by a path of open edges. In terms of $\theta$, $\theta(p)>0$ is equivalent to the existence of an infinite open cluster (with probability one), which naturally leads to a definition of the \emph{critical probability} $p_c = \sup\{ p: \theta (p)=0 \}$.

An argument analogous to that used by Peierls~\cite{Peierls36} to establish the existence of a phase transition in the Ising model can be used to show that $\theta (p)=0$ when $p$ is sufficiently small, and $\theta (p)>0$ when $p$ is sufficiently close to one, which implies that $0<p_c<1$. This fundamental result, proved early on by Broadbent and Hammersley~\cite{BH57,Hammersley57,Hammersley59}, shows that the percolation model undergoes a phase transition and explains the subsequent interest in the subject (see, e.g., \cite{Kesten82,SA94,Grimmett99,BR06}).

When $\theta(p)>0$, it is natural to ask about the multiplicity of infinite open clusters. A simple and elegant proof by Burton and Keane~\cite{BK89} shows that, under quite general conditions, there is almost surely a unique infinite cluster. The same conclusion had however already been reached two years earlier in a joint paper of Chuck's with Michael Aizenman and Harry Kesten~\cite{AKN87}. Their analysis applies to both site and bond models in arbitrary dimension, including long range bond percolation where one considers edges (bonds) between vertices that are not nearest-neighbors. A particularly interesting example of such long-range models, covered in~\cite{AKN87}, is provided by one-dimensional $1/|x-y|^2$ models, where the probability that the bond between vertices $x$ and $y$ is open decays like $1/|x-y|^2$. We'll come back to such models shortly.

An important precursor to \cite{AKN87} is the paper by Chuck and Lawrence Shulman~\cite{NS81} which investigates the number and nature of infinite clusters in a large class of percolation models in general dimension. Using primarily methods from ergodic theory and measure theory, the paper shows that, under general conditions, the number of infinite clusters is either $0$, $1$ or $\infty$. The class of percolation models to which this result applies is characterized by translation invariance, translation ergodicity, and a ``finite energy'' condition which implies that the conditional probability of a local configuration, conditioned on the configuration in the rest of the system, is always strictly between $0$ and $1$. Intuitively, having strictly positive conditional probabilities for local configurations regardless of the rest of the system means that there are no ``prohibited'' (local) configurations. The proof by Burton and Keane~\cite{BK89} also uses this side of the finite energy condition in a crucial way, together with translation invariance and ergodicity. Their proof was extended in~\cite{GKN92}
to more general percolation models, including in particular long-range models. The existing proofs of uniqueness of the infinite cluster may, in fact, be adapted to all ``periodic'' graphs such that the number of vertices within distance $n$ of the origin grows sub-exponentially in $n$. The situation is qualitatively different on trees where, although one still has two phases, one can show that above the critical density $p_c$ there are {\it infinitely many} infinite open clusters instead of just one. Yet another situation arises when one considers the direct product of $\mathbb Z$ and a regular tree, as done in \cite{GN90}, where it is shown that percolation on such a graph has at least {\it three} distinct phases, with the number of infinite clusters being (almost surely) $0$, $\infty$, and $1$, respectively, as the density of open bonds increases. Coauthored by Chuck and Geoffrey Grimmett, \cite{GN90} is the first paper to explore percolation on nonamenable graphs, and the subject has since then attracted much attention. (Roughly speaking, an infinite graph $G$ is nonamenable if, for every finite subset $W$ of $G$, the size of the boundary of $W$ is of the same order as the size of $W$. A typical example of a nonamenable graph is a regular tree graph: as one ``grows'' the tree from the ``root'' adding more vertices to the graph, the number of ``leaves'' is always of the same order as the total number of vertices.) An introduction to percolation on nonamenable graphs and a list of problems (some of which have by now been solved) are contained in~\cite{BS96}, and a somewhat more recent survey can be found in~\cite{Lyons00}.

Another interesting result proved in~\cite{AKN87}, essentially as a corollary of the uniqueness of the infinite open cluster, is the continuity of the percolation function $\theta$, except possibly at the critical density $p_c$. The continuity of the percolation function $\theta$ at $p_c$ is still one of the major open problems in percolation theory. For bond percolation on the square lattice, it was established in a groundbreaking and very influential paper by Harry Kesten~\cite{Kesten80}. Ten years later, Takashi Hara and Gordon Slade~\cite{HS90} proved continuity in dimension $d \geq 19$ (a bound later improved to $d \geq 11$ \cite{FvdH15}) and in more than six dimensions for sufficiently ``spread-out'' models where long-range bonds are allowed. A little later, Barsky, Grimmett and Newman~\cite{BGN91} proved that the probability that there exists an infinite cluster in ${\mathbb N} \times {\mathbb Z}^{d-1}$ is zero for $p=p_c({\mathbb Z}^d)$, but the continuity of $\theta$ at $p_c$ remains a conjecture for percolation models in general dimensions.

The continuity results just described justify the characterization of the percolation phase transition as a {\it continuous} phase transition. Indeed, the percolation phase transition is often considered a prototypical example of a continuous phase transition. However, Chuck and Michael Aizenman~\cite{AN86} showed that the situation is different if one considers one-dimensional $1/|x-y|^2$ models, already mentioned earlier. In such models the percolation function $\theta(p)$ is {\it discontinuous} at $p=p_c$. Chuck and Larry Schulman~\cite{NS86} had already proved that one-dimensional $1/|x-y|^s$ percolation models have a phase transition for all $s \leq 2$, a result that is non-trivial for $s>1$. This is analogous to the occurrence of a phase transition in long-range one-dimensional Ising models with interactions decaying like $1/|x-y|^s$ \cite{Dyson69}. The results of~\cite{AN86} were extended to Ising and Potts models in~\cite{ACCN88}, which provides a proof of a type of discontinuity originally predicted by Thouless~\cite{Thouless69}.

{\it Critical exponents}.
Typically, for statistical mechanical models undergoing a continuous phase transition, close to the critical point, the correlation length and other thermodynamic quantities exhibit power-law behavior in the parameters of the model. The exponents in those power laws are called {\it critical exponents} and their values appear to be largely independent of the microscopic details of the model. Instead, they appear to depend only on global features such as the dimension and symmetries of the model. This phenomenon is called {\it universality} and has a natural explanation within the framework of the {\it renormalization group}, which will be briefly discussed in Section~\ref{subsec:scaling}.

A first quantitative theory of critical phenomena was proposed by Landau~\cite{Landau37}, corresponding to the mean-field approximation that applies to systems on a Bethe lattice~\cite{Bethe35} or in sufficiently high dimensions or with long-range interactions. However, Onsager's exact solution~\cite{Onsager44} of the two-dimensional Ising model and Guggenheim's results~\cite{Guggenheim45} on the coexistence curve of simple fluids showed that critical exponents can take values different from those of Landau's mean-field theory. Indeed, critical exponents should only take their mean-field values above the {\it upper critical dimension} $d_{uc}$, already introduced in Section~I\hspace{-1mm} A. Even today, only a small number of non-mean-field critical exponents have been derived rigorously (including, as we'll discuss later, a couple recently established by Chuck and co-authors in the case of the planar Ising model). But in 1963, a seminal paper by Rushbrooke~\cite{Rushbrooke63} demonstrated how {\it inequalities} between critical exponents can be derived rigorously and exploited fruitfully to study the singularity of thermodynamic functions near the critical point. Since then, numerous such inequalities have been proved for various models of statistical mechanics.

Chuck's contributions to this line of inquiry appeared in several papers~\cite{AN84,Newman86,Newman87,Newman87bis}. Papers~\cite{AN84,Newman86,Newman87} deal with the exponent $\gamma$ for percolation, associated with the expected cluster size $\chi$, namely $\chi(p) \sim (p_c-p)^{-\gamma}$ as $p \nearrow p_c$. In particular, among the results reported in~\cite{Newman86,Newman87} is the fact that in one-dimensional $1/|x-y|^2$ models, where a discontinuous phase transition occurs, $\gamma \geq 2$; conversely, for models such as standard site or bond percolation in dimension $d>2$, where it is believed but not proved that the phase transition is continuous, it is shown that, in order to prove continuity of the percolation function at $p_c$, it would suffice to show that $\gamma<2$. (Incidentally, for percolation in $d=3$, $\gamma$ is numerically estimated to be about 1.7 --- see, e.g., \cite{Stauffer81} and the references given there.) Another result discussed in~\cite{Newman86,Newman87} is the inequality $\gamma \geq 2 (1-1/\delta)$, with $\delta \geq 2$, which improves on the bound $\gamma \geq 1$, proved in~\cite{AN84}.
The critical exponents considered in~\cite{Newman87bis} are the exponent $\beta$, which determines the divergence of the percolation function $\theta$ as $p \searrow p_c$, namely $\theta(p) \sim (p-p_c)^{\beta}$, and the exponent $\delta$, which determines the behavior of the cluster size distribution at $p_c$, namely the probability $P_n(p_c)$ that the cluster of the origin contains exactly $n$ vertices when $p=p_c$. Assuming that the percolation density vanishes at the critical point, the inequality proved in~\cite{Newman87bis} is $\beta \geq 2/\delta$, improving on a previous result by Aizenman and Barsky \cite{AB87} and on $\delta \geq 2$, since $\beta \leq 1$ \cite{CC86}.

Coming back to~\cite{AN84}, besides the proof of the inequality $\delta \geq 1$, already mentioned and soon improved upon, and other results discussed there, the most remarkable and influential contribution of the paper is arguably the introduction of the {\it triangle diagram} $\nabla$ and of the corresponding {\it triangle condition}. The triangle diagram is defined as a sum of two-point functions $\tau$, namely $\nabla = \sum_{x,y} \tau(0,x)\tau(x,y)\tau(y,0)$, and the triangle condition corresponds to the {\it finiteness} of the triangle diagram at $p_c$ (or uniform boundness for $p<p_c$). A main result of~\cite{AN84} is a proof that, in finite-range percolation models, the triangle condition implies $\gamma=1$. The relevance of this result stems from the fact that $1$ is the mean-field value of the exponent $\gamma$. Therefore, if the triangle condition is satisfied for some dimension $d$, this suggests that $d$ is greater than the upper critical dimension $d_{uc}$ and that the mean-field approximation gives the correct prediction for the critical exponents in dimension $d$. Indeed, Barsky and Aizenman \cite{BA91} extended the result of~\cite{AN84} to the critical exponents $\delta$ and $\beta$ by showing that, in percolation models where the triangle condition is satisfied, the exponents $\delta$ and $\beta$ exist and take their mean-field values: $\delta=2$ and $\beta=1$. In particular, the existence of the exponent $\beta$ implies the continuity of the percolation function at $p_c$, i.e., $\theta(p_c)=0$. Another result concerning the triangle condition is the proof by Nguyen~\cite{Ng87} that, if the triangle condition is satisfied, then the {\it gap exponents}, characterizing the divergence of higher moments of the cluster size distribution, assume their mean-field value of $2$.
The triangle condition itself was proved in~\cite{HS89,HS90} for percolation in sufficiently high dimensions for nearest-neighbor models, and above six dimensions for a class of spread-out models. Similar conditions have been subsequently introduced in the literature~\cite{NgY93,BW98,BvdHSS05,BvdHSS05bis}.

{\it First passage percolation}.
\label{subsec:first-passage}
First passage percolation was introduced by Hammersley and Welsh~\cite{HW65} as a percolation-type model with a time dimension that makes it suitable for studying the spread of a disease or a fluid in a porous medium. In the standard version, one assigns to each edge $e$ of ${\mathbb Z}^d$ a nonegative random variable $t(e)$, which is usually interpreted as the passage time of the edge $e$. The {\it passage time} of a path $r$ consisting of edges $e_1,\ldots,e_n$ is $T(r) = \sum_{i=1}^{n} t(e_i)$ and the {\it passage} or {\it travel time} between two vertices $u,v \in {\mathbb Z}^d$ is the infimum of $T(r)$ over all paths $r$ from $u$ to $v$. The stochastic region $\tilde B(t) = \{ x \in {\mathbb Z}^d : T(0,x) \leq t \}$ is the set of vertices that can be reached from the origin by time $t$. The interested reader is referred to~\cite{ADH17} for precise definitions and a comprehensive recent survey.

A main object of interest in first-passage percolation is the set $\tilde B(t)$, and in particular its asymptotic properties when $t$ is large. The interface between $\tilde B(t)$ and its complement will, under natural hypotheses, grow linearly in $t$ with a deterministic shape, while the magnitude of the fluctuations of this interface about its mean shape is believed to be typically of order $t^{\chi}$, with a universal exponent $\chi$ (which should of course depend on the dimension $d$). The study of the fluctuations of growing interfaces is a subject that has attracted considerable attention in the physics literature (see~\cite{KS91} for a review). In the case of first-passage percolation, the first rigorous results were established by Kesten~\cite{Kesten86,Kesten93}, who proved the first rigorous bounds on the variance of the passage time $T(0,x)$. Since then, there has been an extensive literature on both lower and upper bounds for the variance of passage times. A detailed review of that literature can be found in~\cite{ADH17}. In two dimensions, the best lower bound to date was obtained in~\cite{NP95} by Chuck and Marcelo Piza, who showed that the variance of $T(0,x)$ must grow at least as fast as $\log |x|$ and provided the first proof of divergence of the fluctuations of the interface between $\tilde B(t)$ and its complement. Prior to this paper, there had been no proof that the shape fluctuations diverge for any model in any dimension $d>1$.

In related work~\cite{LNP96}, Cristina Licea, Chuck and Marcelo Piza provided the best rigorous results on the {\it wandering exponent} $\xi$, which determines the traverse fluctuations of time-minimizing paths. Like $\chi$, the exponent $\xi$ is expected to depend on the dimension $d$ but to be otherwise universal (e.g., independent of the distribution of the random variables $\tau(e)$). There are, a priori, many possible mathematical definitions of the exponent $\xi$, some based on point-to-plane and some based on point-to-point passage times, but it is believed that they yield the same exponent. Furthermore, it has been conjectured that $\xi(d) \geq 1/2$ for all dimensions $d \geq 2$, with a strict inequality (superdiffusivity) at least for low $d$, and with $\xi(2)=2/3$. In~\cite{LNP96}, Licea, Newman and Piza, working primarily with definitions of $\xi$ of the point-to-plane type, obtained the lower bounds $\xi(d) \geq 1/2$ for all $d$ and $\xi(2) \geq 3/5$. It should be noted that the exponents $\chi$ and $\xi$ had been conjectured in numerous physics papers to satisfy the scaling identity $\chi = 2\xi -1$, irrespective of the dimension (see~\cite{KS91}). This relation was recently proved~\cite{Chatterjee13} by Sourav Chatterjee. Before Chatterjee's paper, the best result was due to Newman and Piza \cite{NP95} who proved that $\chi' \geq 2\xi-1$, where $\chi'$ is an exponent closely related (and perhaps equal) to $\chi$.

Another line of inquiry pursued by Chuck concerns the study of infinite geodesics, motivated by the connection (in $d=2$) between ``bigeodesics'' (i.e., doubly-infinite geodesics) in first-passage percolation and ground states of disordered ferromagnetic spin models~\cite{FLN91,Newman97}. Originating from the physics literature on disordered Ising models is a conjecture that, at least in two dimensions, bigeodesics should not exist. This derives from the conjecture that, in low dimensions (including $d=2$), disordered Ising ferromagnets should have only two (constant) ground states. Indeed, the existence of nonconstant ground states for a disordered Ising model on ${\mathbb Z}^2$ with couplings $J_{e}$ would imply the existence of bigeodesics on a dual square lattice ${{\mathbb Z}^2}^*$ with passage times $\tau(e^*)=J_{e}$, where the dual edge $e^*$ is the perpendicular bisector of $e$.

A celebrated shape theorem~\cite{Richardson73,CD81,Kesten86} states, roughly speaking, that $\tilde B(t)$ behaves like $t B_0 + o(t)$ as $t \to \infty$, where $B_0$ is a convex subset of ${\mathbb R}^2$. Under the assumption that the boundary of $B_0$ has uniformly positive curvature, Cristina Licea and Chuck proved~\cite{LN96} that all infinite geodesics have an asymptotic direction, and that there is a set $D \subset [0,2\pi)$ of full Lebesgue measure such that, for any $\theta \in D$, there are no bigeodesics with one end directed in direction $\theta$. This result provides a partial verification of the conjecture mentioned above; however, a proof that the asymptotic shape $B_0$ has uniformly positive curvature seems to be still out of reach.

In~\cite{HN97}, Doug Howard and Chuck introduced a Euclidean version of first-passage percolation where the vertex set of the lattice ${\mathbb Z}^d$ is replaced by the set of points of a homogeneous Poisson point process on ${\mathbb R}^d$. For this model, in~\cite{HN97} they proved a shape theorem, with asymptotic shape given by a Euclidean ball (due to the isotropy of the model), and the nonexistence of certain (geometrically non-realistic) doubly infinite geodesics. In~\cite{HN99}, they continued their study of this new model, showing in particular that, for any dimension, with probability one, every semi-infinite geodesic has an asymptotic direction and every direction has at least one semi-infinite geodesic (starting from each Poisson point). In~\cite{HN01}, they proved the bounds $\chi \leq 1/2$ and $\xi \leq 3/4$ for Euclidean first-passage percolation in all dimensions, together with other interesting results concerning spanning trees of semi-infinite geodesics and related random surfaces.

Chuck's work on first passage percolation has been and continues to be very influential. His results on the subject continue to be cited, and the ideas and methods developed by Chuck and coauthors have been used and extended in numerous papers (for a small selection, see~\cite{FP05,Hoffman08,GM10,CP11,DH13,AD13,DH17,BSS17+}).

\subsection{Limit Theorems and Continuum Scaling Limits}
\label{subsec:scaling}

Over the past two decades, besides working on disordered systems and nonequilibrium dynamics, Chuck has dedicated a lot of effort to the study of {\it continuum scaling limits}, particularly in the case of spanning trees, percolation and the Ising model. Chuck's interest in scaling limits, and in particular in the Euclidean fields emerging from such limits, has deep roots, going all the way back to his Ph.D. dissertation on Quantum Field Theory \cite{Newman72} (see also~\cite{Newman73}).

A main goal of both probability theory and statistical physics is to understand and describe the behavior of random systems with a very large number of ``degrees of freedom.'' In field theory, one deals with an infinite number of degrees of freedom, and there is indeed a deep connection between the theory of critical phenomena and field theory.
This connection is particularly salient in the renormalization group approach to these theories in which, broadly speaking, one analyses the behavior of specific observable quantities at different scales (see, e.g., \cite{WK74,Wilson75,Wilson83}). As more scales are included, the behavior of these observable quantities is random, and the quantities themselves need to be ``renormalized'' to ``tame'' their fluctuations. \emph{Renormalizable} systems are such that this renormalization procedure is possible and one can take a limit over all the scales of the system. In field theory, a momentum cutoff introduced to reduce the number of degrees of freedom is sent to infinity as higher and higher energies are taken into consideration. In statistical mechanics, a version of the renormalization approach can be implemented with a continuum scaling limit in which some elementary scale of the system (e.g., the lattice spacing in lattice-based models such as percolation and the Ising model) is sent to zero. Renormalizable systems may possess a certain amount of \emph{scale invariance} leading, in the limit of large scales or large momenta, to models that are (statistically) \emph{self-similar}.

Thanks to the work of Polyakov~\cite{Polyakov70} and others~\cite{BPZ84,BPZ84bis}, it was understood by physicists since the early seventies that critical statistical mechanical
models should typically possess continuum scaling limits with a global conformal invariance that goes beyond simple scale invariance. Following this observation, in two dimensions, conformal methods were applied extensively to Ising and Potts models, Brownian motion, the self-avoiding walk, percolation, diffusion limited aggregation, and more. The large body of knowledge and techniques (mainly non-rigorous) that resulted from these efforts goes under the name of Conformal Field Theory (CFT) --- see~\cite{DFMS97} for an extensive treatment.

{\it Normal fluctuations, the Gaussian fixed point, and beyond}.
The renormalization group idea of gradually including the effect of more scales or degrees of freedom has a parallel in the theory of limit theorems in probability theory (see \cite{BS73,J-L75,J-L76,J-L01}).

In two influential papers~\cite{Newman80,Newman83}, Chuck provided general conditions for the central limit theorem to apply to positively correlated random variables. Normal fluctuations are expected to occur in statistical mechanical systems away from a critical point, in which case, observables such as the renormalized magnetization in Ising ferromagnets are expected to converge to Gaussian random variables. In renormalization group terminology, the system converges to the {\it Gaussian fixed point}. Chuck's results have several important applications, including to magnetization and energy fluctuations in a large class of Ising ferromagnets, infinite cluster volume and surface density fluctuations in percolation models, and boson field fluctuations in (Euclidean) Yukawa quantum field theory models.

While Gaussian limits are very common, non-Gaussian limits can arise when dealing with random variables that are strongly positively correlated. This is believed to happen in various models studied in statistical mechanics, in particular in (ferromagnetic) Ising systems at the critical point. In the 1970s (extending earlier results of Simon and Griffiths \cite{SG73}), Chuck and co-authors studied the emergence of non-Gaussian limits in Curie-Weiss (mean-field) models~\cite{EN78,EN78bis,ENR80}. Part of the motivations for choosing a relatively simple class of systems such as mean-field models came from the fact that, at the time, the existence of non-Gaussian limits was an open problem, and seemed out of reach, for almost all non-trivial Ising systems. The situation has changed since then, and Chuck has contributed significantly to this change. As we will discuss at the end of this section, more than thirty years after his work on Curie-Weiss models, Chuck and a different set of co-authors were able to provide substantial contributions to the study of other non-Gaussian limits, this time in the context of the critical planar Ising model. Meanwhile, new interesting results on other mean-field models, some of which are analogous to those of~\cite{EN78,EN78bis}, have also emerged --- see the review article 
by Alberici, Contucci and Mingione in this volume, where the authors describe their recent work on mean-field monomer-dimer models.

{\it Critical and near-critical percolation}.
The renormalization group approach has greatly improved our understanding of critical phenomena, but from a mathematical perspective it remains to this day largely non-rigorous.
In the late 1990s, Chuck contributed to the development of a new framework~\cite{Aizenman98,AB99,ABNW99}. That framework and, in particular, the introduction of the Shramm-Loewner Evolution (SLE) by Oded Schramm~\cite{Schramm00} provided a new, mathematically rigorous, approach to study the geometry of critical systems on the plane. This new approach consists in viewing cluster boundaries in models such as Ising, Potts and percolation models, or loops in the $O(n)$ model, as random interfaces with a distribution that depends on the scale of the system under consideration (i.e., the lattice spacing), and in analyzing the continuum limit as the scale of the system is sent to zero. Schramm realized~\cite{Schramm00} that, at criticality, these interfaces become, in the continuum limit, random planar curves whose distributions can be uniquely identified thanks to their conformal invariance and a certain ``Markovian" property. There is a one-parameter family of SLEs indexed by a positive real number $\kappa$: SLE$_\kappa$ is a random growth process based on a Loewner chain driven by a one-dimensional Brownian motion with speed $\kappa$ \cite{Schramm00,RS05}. Lawler, Schramm and Werner \cite{LSW01,LSW01bis,LSW02,LSW04}, and Smirnov and Werner \cite{SW01}, among others, used SLE to confirm several results that had appeared in the physics literature, and to prove new ones.

In particular, substantial progress was made, thanks to SLE, in understanding the fractal and conformally invariant nature of (the scaling limit of) large percolation clusters, starting with the work of Schramm~\cite{Schramm00} and Smirnov~\cite{Smirnov01}, which identified the scaling limit of critical percolation interfaces with SLE$_6$. One of Chuck's early contributions to this area of research, in collaboration with one of us (FC), was the first complete and self-contained proof of the convergence of a percolation interface in critical site percolation on the triangular lattice to SLE$_6$ \cite{CN07}. The proof, which follows roughly Smirnov's strategy~\cite{Smirnov01} with some important differences, appeared first in an appendix of an unpublished preprint~\cite{CN05}. It was later deemed worth publishing as a separate paper~\cite{CN07} because of the importance of the result and of its many applications, which include the rigorous derivation of various percolation critical exponents~\cite{SW01}, of Watt's crossing formula~\cite{Dubedat06}, and of Schramm's percolation formula~\cite{Schramm01}.

The convergence result proved in~\cite{CN07} is also a crucial ingredient in the derivation of the {\it full scaling limit} of critical percolation~\cite{CN04,CN06}, obtained by FC and Chuck by taking the continuum limit of the collection of all percolation interfaces at once. The resulting object, called the {\em continuum nonsimple loop process} in~\cite{CN04,CN06}, is invariant under scale and conformal transformations, and inherits a spatial Markov property from the discrete percolation process. As shown in~\cite{CN08}, the continuum nonsimple loop process constructed in~\cite{CN04,CN06} is a {\it Conformal Loop Ensemble} (CLE) \cite{Sheffield09,SW12}. The construction of the continuum scaling limit of planar critical percolation~\cite{CN04,CN06} has had several interesting applications, including~\cite{BN11,Yao14,Yao18,Yao18+,Yao18++}.

Another notable application of~\cite{CN04,CN06} is in the construction of the {\it near-critical} scaling limit of planar percolation in which the percolation density $p$ approaches the critical value $p_c$ ($p_c=1/2$ for site percolation on the triangular lattice) according to an appropriate power of the lattice spacing $a$, $p=p_c+\lambda a^{\alpha}$, as $a \searrow 0$. With an appropriate choice of $\alpha$ ($\alpha=3/4$ for site percolation on the triangular lattice), one can show that the scaling limit leads to a one-parameter family, indexed by $\lambda$, of limits that are not scale invariant (except in the critical case, $\lambda=0$) but retain some of the properties of the critical scaling limit (see, e.g., \cite{CJM09}). In \cite{CFN06,CFN06bis}, Camia, Fontes and Newman proposed an approach to construct the one-parameter family of near-critical scaling limits of planar percolation based on the critical full scaling limit and the ``Poissonian marking'' of some special (``macroscopically pivotal'') points. This approach leads to a conceptual framework that can describe not only the scaling limit of near-critical percolation but also of related two-dimensional models such as dynamical percolation, the minimal spanning tree and invasion percolation. Partly inspired by~\cite{CFN06,CFN06bis}, Garban, Pete and Schramm later constructed those scaling limits in~\cite{GPS13}. Although Garban, Pete and Schramm use a notion of full scaling limit different from that of~\cite{CN04,CN06}, they do take advantage of some of the results proved in~\cite{CN06}, in particular to resolve issues of uniqueness, measurability, and conformal invariance of the limit.

Besides the applications mentioned above, a remarkable aspect of~\cite{CN04,CN06} is that they contain at once the first example of a nested CLE, including a description of some of its crucial properties, as well as the first rigorous construction of a CLE as a scaling limit, before the concept of CLE was formalized in full generality in~\cite{Sheffield09}. Conformal Loop Ensembles have been extensively studied for their intrinsic interest (see, e.g., \cite{SW12,SSW09,WW13,MSW14,MWW15,MWW16,KW16,MSW17,MW17+}) and for their applications to CFT (see, e.g., \cite{Doyon12,Doyon13,Doyon14}), and we don't attempt to provide here a comprehensive list of references.

{\it Critical and near-critical Ising model}. 
The Ising model~\cite{Ising25}, suggested by Lenz~\cite{Lenz20} and cast in its current form by Peierls~\cite{Peierls36}, is one of the most studied models of statistical mechanics. Its two-dimensional version has played a special role since Peierls' proof of a phase transition~\cite{Peierls36} and Onsager's calculation of the free energy~\cite{Onsager44}. This phase transition has become a prototype for developing new techniques. Its analysis has helped test a fundamental tenet of critical phenomena, that near-critical physical systems are characterized by a \emph{correlation length}, which provides the natural length scale for the system and diverges when the critical point is approached.

Substantial progress in the rigorous analysis of the two-dimensional Ising model at criticality was made by Smirnov~\cite{Smirnov10} with the introduction and scaling limit analysis of \emph{discrete holomorphic observables}. These have proved extremely useful in studying the Ising model in finite geometries with boundary conditions and in establishing conformal invariance of the scaling limit of various quantities, including the energy density~\cite{Hongler10,HS13} and spin correlation functions~\cite{CHI15}. (An independent derivation of critical Ising correlation functions in the plane was obtained in~\cite{Dubedat11}.)

In~\cite{CN09}, Chuck and one of us (FC) proposed a strategy to obtain the continuum scaling limit of the renormalized Ising magnetization field using the FK random cluster representation of the Ising model (see, e.g., \cite{Grimmett06}) 
and a new tool called, in later papers, a {\it conformal measure ensemble}. The strategy, which involves coupling conformal measure and conformal loop ensembles and leads to a geometric representation of the Ising magnetization in the continuum reminiscent of the Edwards-Sokal coupling~\cite{ES88} was not fully carried out in~\cite{CN09}. (In particular, no conformal  measure ensemble is constructed in~\cite{CN09}.) Nevertheless, the paper was the starting point for the subsequent work of Camia, Garban, Newman and of Camia, Jiang, Newman, which we now discuss.

In~\cite{CGN15} FC, Christophe Garban and Chuck, partly following the strategy of~\cite{CN09} but without the use of conformal measure ensembles, showed that, in two dimensions at the critical point, when properly renormalized, the Ising magnetization field has a continuum scaling limit $\upPhi^0$ which satisfies the expected properties of conformal covariance. $\upPhi^0$ is a (generalized, Euclidean) random field on $\mathbb{R}^2$ --- i.e., for suitable test functions $f$ on $\mathbb{R}^2$, one can construct random variables $\upPhi^0(f)$, formally written as $\int_{\mathbb{R}^2}\upPhi^0(x)f(x)\upD x$. The tail behavior of the field $\upPhi^0$, obtained in~\cite{CGN16}, shows that $\upPhi^0$ is not a Gaussian field. (This follows also from the behavior of its correlations, which do not obey Wick's theorem --- see~\cite{CHI15}). Another significant contribution of~\cite{CGN16} is the construction of the continuum scaling limit of the magnetization field for the near-critical Ising model with external magnetic field $ha^{15/8}$ on the rescaled lattice $a\mathbb{Z}^2$ as $a \searrow 0$. As stated in~\cite{CGN15}, it was expected that the truncated correlations of the resulting field $\upPhi^h$ would decay exponentially whenever $h \neq 0$. A proof of that statement is provided in~\cite{CJN18+} together with a rigorous proof that the critical exponent for how the correlation length diverges as $h \searrow 0$ is $8/15$, and related scaling properties of $\upPhi^h$. The related magnetization critical exponent $\delta$, which determines how, at the critical temperature, the magnetization depends on the external magnetic field, was rigorously shown to be equal to 15 in a joint paper by FC, Garban and Chuck~\cite{CGN14}.

It is worth pointing out that, while the concept of conformal measure ensemble is not used in~\cite{CGN15,CGN16}, it does play a crucial role in~\cite{CJN18+}.
Indeed, a surprising contribution of~\cite{CJN18+} is the demonstration that conformal measure ensembles coupled with the corresponding conformal loop ensembles can be useful in analyzing near-critical scaling limits. The construction of conformal measure ensembles and their coupling to CLE$_{\kappa}$ is carried out for $\kappa=6$ and $16/3$ in the article by FC, Ren\'e Conijn and Demeter Kiss in Volume 2 of this Festschrift, 
where the authors also obtain the geometric representation of the Ising magnetization in the continuum mentioned earlier. We note that, as pointed out in \cite{CN09}, in addition to their utility for critical and near-critical two-dimensional models, measure ensembles may be more extendable than loop ensembles to scaling limits in dimensions $d>2$.

It seems fitting to conclude this introduction with some comments on the relevance of~\cite{CGN15,CGN16}
and \cite{CJN18+} for (constructive) quantum field theory, one of Chuck's early scientific interests.
Euclidean random fields, such as $\upPhi^h$, on the Euclidean ``space-time'' $\mathbb{R}^d:=\{x=(x_0,w_1,\ldots,w_{d-1})\}$ are related to quantum fields on relativistic space-time, $\{(t,w_1,\ldots,w_{d-1})\}$, essentially by replacing $x_0$ with a complex variable and analytically continuing from the purely real $x_0$ to a pure imaginary $(-it)$ --- see~\cite{OS73,GJ87}. In this context, non-Gaussian Euclidean fields such as those discussed in~\cite{CGN15,CGN16} are of particular interest since Gaussian Euclidean fields correspond to non-interacting (and therefore trivial) quantum fields --- see, e.g., \cite{FFS92}. The construction of interacting field theories has been from the start one of the main goals of constructive field theory, albeit in the physically interesting case of four dimensions rather than in the two-dimensional setting of~\cite{CGN15,CGN16}. Notwithstanding, one major reason for interest in $\upPhi^h$ is that the associated quantum field theory was predicted by Zamolodchikov~\cite{Zamolodchikov89} to have remarkable properties including a ``particle content'' of eight particles whose masses are related to the exceptional Lie algebra $E_8$ --- see \cite{Delfino04,BG11,MM12}. One of the main contributions of~\cite{CJN18+}
is to prove a strictly positive lower bound on all masses (i.e., a \emph{mass gap}) predicted by Zamolodchikov's theory. This is a first natural step in the direction of a rigorous analysis of Zamolodchikov's theory, to which we hope Chuck will provide further insightful contributions.

The description of Chuck's contributions in this introduction is unavoidably biased and incomplete, and we apologize to all of his co-authors whose work has not been mentioned. The responsibility rests with Chuck for being so unreasonably productive.

%
%

\end{document}